\documentclass{article}
\usepackage{amsmath,amssymb,latexsym}

\begin{document}
\title{A solution to a problem of Fermat, on two numbers of which the sum
is a square and the sum of their squares is a biquadrate, inspired by the
Illustrious La Grange\footnote{Delivered to the
St.--Petersburg Academy June 5, 1780.
Originally published as
{\em Solutio problematis Fermatiani de duobus numeris, quorum summa sit quadratum, quadratorum vero summa biquadratum, ad mentem illustris La Grange adornata},
M\'emoires de l'Acad\'emie Imp\'eriale des Sciences de St.-P\'etersbourg \textbf{10} (1826),
3--6, and
republished in \emph{Leonhard Euler, Opera Omnia}, Series 1:
Opera mathematica,
Volume 5, Birkh\"auser, 1992. A copy of the original text is available
electronically at the Euler Archive, at http://www.eulerarchive.org. This paper
is E769 in the Enestr\"om index.}}
\author{Leonhard Euler\footnote{Date of translation: April 23, 2006.
Translated from the Latin
by Jordan Bell, 4th year undergraduate in Honours Mathematics, School of Mathematics and Statistics, Carleton University,
Ottawa, Ontario, Canada.
Email: jbell3@connect.carleton.ca.}}
\date{}
\maketitle

\S 1. For the solutions of this problem, which so far have been
conveyed in public, the Illustrious La Grange rejects them most deservedly,
because they are taken altogether by chance and roving efforts,
so that we are not able to be certain of all solutions. 
I am confident for the following analysis to give satisfaction
to this desire.

\S 2. Were $x$ and $y$ two numbers that are being sought, such that it mught
be that
$x+y=\square$ and $xx+yy=\square^2$, if for the first condition we take
$x=pp-qq$ and $y=2pq$, then it will be $xx+yy=(pp+qq)^2$. Then if
in turn it were set $p=rr-ss$ and $q=2rs$, it would be
$pp+qq=(rr+ss)^2$, and thus $xx+yy=(rr+ss)^4$, as is required. Then moreover
it will be $x=r^4-6rrss+s^4$ and $y=4rs(rr-ss)$.

\S 3.For the first condition therefore the sum of the numbers will be
\[
x+y=r^4+4r^3s-6rrss-4rs^3+s^4,
\]
which formula is therefore to be made into a square. To this end, 
as nothing else presents itself for this effort,  I represent this expression below in this form:
\[
x+y=(rr+2rs-ss)^2-8rrss,
\]
such that now this form: $AA-2BB$ ought to be turned into a square,
where it is to be taken $A=tt+2uu$ and $B=2tu$; for then it would be
\[
AA-2BB=(tt-2uu)^2.
\]

\S 4. Now in place of $A$ and $B$ we may write our values and we will have
$rr+2rs-ss=tt+2uu$ and $2rs=2tu$, and so in this way the sum
of our numbers will be $x+y=(tt-2uu)^2$, and therefore now both conditions
will be satisfied,
providing that the formulas found above are satisfied.

\S 5. Seeing moreover that the two products $rs$ and $tu$ must be equal to
each other, in place of the letter $s$, unity may be taken without risk.
However, now for $r$ fractions may be produced, but this does
not impede the solution at all, because a solution found in a fraction
can easily be reduced to integers. It will therefore be $r=tu$, whose
value substituted into the other equation gives $ttuu+2tu-1=tt+2uu$,
and thus the problem is reduced to finding the proper relation between
$t$ and $u$. We wish therefore to define either $t$ by $u$ or $u$ by $t$;
the resolution of this equation by quadrature supplies the following two
formulaes:
\[
t=\frac{u \pm \sqrt{2u^4-1}}{1-uu} \quad \textrm{and} \quad u=\frac{t \pm
\sqrt{t^4-2}}{2-tt}.
\]
Moreover then, at once values of the radicals are produced for the following
use, so that we will not require the extraction of roots beyond this.
From the first it will be $\sqrt{2u^4-1}=t(1-uu)-u$; from
the other $\sqrt{t^4-2}=u(2-tt)-t$. The benefit of this method
is that from both formulas pairs of values are produced.

\S 6. We shall begin with the first formula, because the case $u=1$ at once
comes to the eyes. Seeing indeed in this case that the denominator 
$1-u$ vanishes, a well-known remedy is to be turned to, in which
it ought to be put $u=1-\omega$, where $\omega$ denotes a vanishing quantity,
so that higher powers of it can be safely rejected. Then it will therefore
be $2u^4=2-8\omega$, and thus $\sqrt{2u^4-1}=\sqrt{1-8\omega}=1-4\omega$ and
$1-uu=2\omega$, and hence it is gathered $t=\frac{3}{2}$, whose
value substituted into the other formula gives
$\sqrt{2t^4-2}=\frac{7}{4}$.

\S 7. The other equation is now progressed to, for which we now have the
values $u=1$ and $t=\frac{3}{2}$, and because a pair of values is held,
we may elicit a new value for $u$, namely $u=-13$. We may bring
this value to the first formula, for which we now know the value to
be $t=\frac{3}{2}$. Because it is known that
\[
\sqrt{2u^4-1}=t(1-uu)-u,
\]
as it is $u=-13$ and $t=\frac{3}{2}$, from this it will be
$\sqrt{2u^4-1}=-239$. Now indeed the above equation supplies us with a new
value for $t$, namely $t=-\frac{113}{84}$.

\S 8. In a similar way  we may infer the value in the other equation,
and because it was $u=-239$, we may deduce
\[
\sqrt{t^4-2}=u(2-tt)-t=-\frac{311485}{7056},
\]
whose value applied to the other radical gives us a new value for $u$, namely
$u=\frac{301993}{1343}$. And if again this value is used in the first
formula, again a new value for $t$ is secured,  
and thus it is permitted to procede as far as we please. Soon however,
because of the immense numbers, the labour is compelled to stop.

\S 9. The strength of this new method consists therefore in that
for each value of $t$ a pair of values of $u$ correspond, and in the same
way for each of $u$ two a pair of values of $t$; therefore, with us having
carried
this out, we may exhibit this for 
inspection:
\begin{eqnarray*}
u=1;&t=\frac{3}{2},\\
u=-13,&t=-\frac{113}{84},\\
u=\frac{301993}{1343},&
\end{eqnarray*}
of which the values are able to be combined with the two adjacent in any way.
From this moreover the two values of the numbers $x$ and $y$ that are being
searched for may be determined in this way
\begin{align*}
x=t^4u^4-6ttuu+1\\
y=4tu(ttuu-1).
\end{align*}
It is also easily seen that in this way all possible solutions ought necessarily
to be produced.

\S 10. This is especially noteworthy because the values found successively 
for the letters
$t$ and $u$ progress by a singular order, so that from each the following are
able to be defined easily. For whatever values are had for $t$ and $u$
which satisfy the formula $t=\frac{u \pm \sqrt{2u^4-1}}{1-uu}$, with
it $\sqrt{2u^4-1}=t(1-uu)-u$, from the varying radical above the
other value for $t$ may be elicited, which if we put equal to $t'$,
it will
then be $t'(1-uu)=2u-t(1-uu)$, and thus $t'=\frac{2u}{1-uu}-t$.

\S 11. In the same way for, when the values of $t$ and $u$ are known,
by the other formula $u=\frac{t \pm \sqrt{t^4-2}}{2-tt}$, because
$\sqrt{t^4-2}=u(2-tt)-u$, the other value for $u$ will be able to be elicited,
which if it is put equal to $u'$, it will be
\[
u'(2-tt)=2t-u(2-tt), \quad \textrm{and thus} \quad u'=\frac{2t}{2-tt}-u.
\]
With these values known, by the other formula again
other new values will be able to be found, which if they
are designated successively by $t'',u''; t''',u'''$; etc. by
$t'=\frac{2u}{1-uu}-t$ and $u'=\frac{2t}{2-tt}-u$, and in a similar way we will
have 
$t''=\frac{2u'}{1-u'u'}-t'$ and $u''=\frac{2t'}{2-t't'}-u'$, then indeed
$t'''=\frac{2u''}{1-u''u''}-t''$ and $u'''=\frac{2t''}{2-t''t''}-u''$; 
and thusly hereafter.

\end{document}